\magnification\magstep1
\overfullrule = 0pt

\def\n{\noindent}
\def\qed{{\hfill{\vrule height7pt width7pt
depth0pt}\par\bigskip}} 
\def\ie{{\it i.e. } }

 \def\qed{{\hfill{\vrule height7pt width7pt
depth0pt}\par\bigskip}}

\def\eg{ {\it e.g.\ }}

\def\calk{{\cal K}}
\def\K{\calk}

\def\NN{{\rm I}\!{\rm N}}
\def\CC{ \;æ{}^{ {}_\vert }\!\!\!{\rm C}}

\centerline{\bf Similarity problems and length}
\bigskip
\centerline{ by Gilles Pisier\footnote*
{Supported in part by the NSF   and
 the Texas Advanced Research
 Program 010366-163. This is an expanded
version of a
lecture given as keynote speaker at the 
January 2000 ICMAA in Taiwan.}}
\centerline{ Texas A\&M University}
\centerline{College Station, TX 77843, U. S. A.}
\centerline{and}
\centerline{Universit\'e Paris VI}
\centerline{Equipe d'Analyse, Case 186, 75252}
\centerline{ Paris Cedex 05, France}\bigskip

\bigskip\bigskip
{\narrower\narrower\medskip\noindent {\bf Abstract.} 
This is a survey of the author's recent results on the Kadison and Halmos
similarity problems and the closely connected notion of ``length'' of an
operator algebra. \bigskip} 

\baselineskip = 18pt

We start by recalling a well known conjecture formulated by Kadison [Ka]
in 1955.

\proclaim Kadison's similarity problem. Let $A$ be a unital $C^*$-algebra
and let $u\colon \ A\to B(H)$ ($H$ Hilbert) 
be a unital homomorphism
(\ie \ we have $u(1)=1$ and $u(ab) = 
u(a) u(b)$ $\forall~a,b\in A$). Show that if
$u$ is bounded, then $u$ is similar to a $*$-homomorphism, 
\ie 
$\exists~\xi\colon \ H\to H$ invertible such that $u_\xi\colon \ a\to
\xi^{-1}u(a)\xi$ is a $*$-homomorphism ($= C^*$-representation).

Explicitly, the conclusion means that
$$\xi^{-1}u(a^*)\xi =
(\xi^{-1}u(a)\xi)^*,\leqno
\forall~a\in A$$ when this holds, Kadison
calls $u$ ``orthogonalizable''. Many partial
results are known, mainly due to Erik
Christensen ([C1--C4]) and Uffe Haagerup
([H1]). In particular, they established (see
[C3] and [H1]) this conjecture for {\it
cyclic\/} homomorphisms, \ie when $u$ admits
a cyclic vector $h$ in $H$ (= a vector $h$
such that $\overline{u(A)h} = H$) or more
generally when $u$ admits a finite cyclic set
$h_1,\ldots, h_n$ (so that we have
$\overline{u(A)h_1 +\cdots+ u(A)h_n} = H$). 

In addition, the Kadison conjecture is known in the following cases:
\item{(i)} $A$ is commutative.
\item{(ii)} $A$ is the unitization,
 denoted by $\widetilde \K$, of the
$C^*$-algebra, denoted by $\K$, of all compact
operators on
$\ell_2$, or more generally when
$A$ is nuclear (see [C2]).
\item{(iii)} $A = B({\cal H})$  or
more generally when $A$ has no tracial states
(see [H1]).
\item{(iv)} $A = \widetilde \K \otimes B$ with
$B$ arbitrary unital
$C^*$-algebra.
\item{(v)} $A$ is a $II_1$-factor with Murray and von~Neumann's property
$\Gamma$ (see [C4]) for instance when $A$ is the so-called hyperfinite
$II_1$-factor (= infinite tensor product of $2\times 2$ matrices with
normalized trace).\medskip

\n In sharp contrast, the conjecture is still open when $A$ is the
reduced $C^*$-algebra of the free group with $\ge 2$ generators, or even
when
$$A = \left(\oplus \sum_{n\ge 1}M_n\right)_\infty = \{x = (x_n)\mid x_n\in
M_n \ \ \sup \|x_n\| < \infty\}.$$
Kadison formulated his conjecture as the $C^*$-algebraic version of a well
known problem (at the time of his writing): \ are all uniformly bounded
group representations similar to unitary representations (= unitarizable).
While a counterexample to that question was soon 
found ([EM], see also [P7] for more recent
results on this theme),  Kadison's conjecture
remained open. Recently, it became entirely
clear that his conjecture is equivalent to
another important open question, the derivation
problem, itself a crucial problem in the
cohomology theory of operator algebras (cf.\
[SS]).

\proclaim Derivation problem. Let $\pi\colon \ A\to B(H)$ be a
$*$-homomorphism (= representation) on a $C^*$-algebra $A$. Let
$\delta\colon \ A\to B(H)$ be a $\pi$-derivation (\ie 
$\delta(ab) =
\pi(a) \delta(b) + \delta(a) \pi(b)$). Show that the boundedness of
$\delta$ (which is actually automatic here) implies that $\delta$ is
inner, which means: \ $\exists~T \in B(H)$ such that $\delta(a) = \pi(a)T
- T\pi(a)$ $\forall~a\in A$. 
We set $$\delta_T(a) = \pi(a)T-T\pi(a).$$

The connection between the two problems is simple. Intuitively derivations
appear as ``infinitesimal generators'' for homomorphisms. More
elementarily, if $\delta$ is as above then
$$u(a) = \left(\matrix{\pi(a)&\delta(a)\cr 0&\pi(a)\cr}\right)$$
is a homomorphism into $B(H\oplus H) = M_2(B(H))$.
  Kirchberg [Ki] recently proved that the $C^*$-algebras which satisfy
the derivation problem are exactly the same as
those which satisfy Kadison's conjecture, but
it is still open whether this class is that of
all
$C^*$-algebras!

We now turn to a key notion to study these problems: \ ``complete
boundedness'' (see [Pa1]).

\proclaim Definition. Let $E\subset B(H)$ and
$F\subset B(K)$ be operator spaces, consider a
map
$$\matrix{B(H)&&B(K)\cr
\cup&&\cup\cr
E&{\buildrel u\over \longrightarrow}&F\cr}$$
For any $n\ge 1$, let
$M_n(E) = \left\{(x_{ij})_{ij\le n}\mid
x_{ij}\in E\right\}$ be the space of $n\times
n$ matrices with entries in $E$. In particular
we have a natural identification
$M_n(B(H)) \simeq B(\ell^n_2(H))$
where $\ell^n_2(H)$ means $\underbrace{H\oplus H\oplus \cdots \oplus H}_{n
\ {\rm times}}$. Thus, we may equip $M_n(B(H))$ and a fortiori its subspace
$M_n(E) \subset M_n(B(H))$
with the norm induced by
$B(\ell^n_2(H)).$
Then, for any $n\ge 1$, the linear map $u\colon \ E\to F$ allows to define
a linear map
$u_n\colon \ M_n(E)\longrightarrow
M_n(F)$   by setting
$$u_n \left(\matrix{&\vdots\cr\ldots&x_{ij}&\ldots\cr &\vdots\cr}\right) =
\left(\matrix{&\vdots\cr \ldots&u(x_{ij})&\ldots\cr &\vdots\cr}\right).$$
A map $u\colon\ E\to F$ is called completely bounded (in short c.b.) if
$$\sup_{n\ge 1} \|u_n\|_{M_n(E)\to M_n(F)}<\infty.$$
We define
$$\|u\|_{cb} = \sup_{n\ge 1} \|u_n\|_{M_n(E)\to M_n(F)}$$
and we denote by
$cb(E,F)$
the Banach space of all c.b.\ maps from $E$ into $F$ equipped with the
c.b.\ norm.

This concept is fundamental in the currently
very actively developed theory of operator
spaces, see [P8].

\proclaim Theorem 1 (Haagerup 1983, [H1]). In the situation of Kadison's
similarity problem, $u$ is similar to a $*$-homomorphism iff $u$ is c.b.
Moreover we have  $$\|u\|_{cb} =
\inf\{\|\xi^{-1}\| \ \|\xi\|\mid \ u_\xi\ 
*-{\rm homomorphism}\}.$$

For derivations, the analogous result is the
following.

\proclaim Theorem 2 (Christensen 1977, [C5]).
In the derivation problem,
$\delta$ is inner iff $\delta$ is c.b. Moreover, we have
$$\|\delta\|_{cb} = \inf\{2\|T\|\mid \delta
= \delta_T\}.$$

Vern Paulsen generalized
Haagerup's result to the non-self-adjoint case:

\proclaim Theorem 3 (Paulsen 1984, [Pa2]). Let $A$ be a unital operator
algebra (\ie  we assume only that $A$ is a
closed subalgebra of $B({\cal H})$ with $I\in
A$). Consider again a homomorphism 
$u\colon\ A\to B(H)$. Then
$\|u\|_{cb} <
\infty$ iff $u$ is similar to a completely
contractive homomorphism, \ie 
$\exists~\xi\colon \ H\to H$ invertible such
that $u_\xi\colon \ a\to \xi^{-1}u(a)\xi$
satisfies
$\|u_\xi\|_{cb} = 1$. Moreover, we have
$$\|u\|_{cb} = \inf\{\|\xi\| \ \|\xi^{-1}\|\mid \|u_\xi\|_{cb} = 1\},$$
and this infimum is attained.

It is easy to see that if $A$ is a $C^*$-algebra then
 $$(\|u\|_{cb}=1)
\Leftrightarrow (\|u\|=1) \Leftrightarrow 
(u \ \hbox{is\  a\  $*$-homomorphism}).$$
 This
explains why Theorem~3 contains Theorem~1 as a
special case. The preceding result leads us
naturally to enlarge our investigation to the
non-self-adjoint case as follows.
\medskip

\n {\bf Generalized Similarity Problem.} Which unital operator algebras
$A$ have the following property denoted by
(SP) ?
\item{(SP)} Any bounded homomorphism $u\colon \
A\to B(H)$ ($H$ arbitrary Hilbert space) is
c.b.

\n Loosely speaking, this property (SP)
 could be described as ``automatic
complete boundedness'' in analogy 
with the field of automatic continuity
for homomorphisms between 
Banach algebras (see [DW]).
\medskip

\n {\bf Example.} The most natural example of a non-self-adjoint algebra is
the disc algebra $A=A(D)$ which can be
described as the completion of the set of all
polynomials $P$ for the norm
$$\|P\|_\infty = \sup\{|P(z)|\mid z\in \partial D\}.$$
We consider $A(D)$ as an operator subalgebra of
the commutative
$C^*$-algebra $C(\partial D)$.

\n Consider a {\it fixed\/} operator $x\in
B(H)$. Let
$$u^x\colon \ P\to P(x) \in B(H)$$
be the homomorphism of evaluation at this fixed $x$. Then $u^x$ is bounded
iff $x$ is polynomially bounded, \ie  $\exists~C$ such that
$$\|P(x)\|\le C\|P\|_\infty.\leqno (1)\qquad\qquad\forall~P$$
On the other hand, it follows from Paulsen's similarity criterion
(Theorem~3 above) that $u^x$ is c.b.\ iff $x$ is similar to a contraction
(which means $\exists~\xi\colon \ H\to H$ invertible such that
$\|\xi^{-1}x\xi\| \le 1$). Indeed, when $\|x\|\le 1$, von~Neumann's
classical inequality shows that (1) holds and actually also (Sz.-Nagy's
dilation) that $\|u^x\|_{cb}=1$. Thus it is the same to ask whether $A(D)$
satisfies (SP) or to ask whether
 any polynomially bounded  operator $x$ is
similar to a contraction. This was a well known problem originally
formulated by Halmos in a landmark 1970 paper
[Ha]. We have recently given a counterexample
as follows.

\proclaim Theorem 4 (1997, [P1]). For any $c>1$ there is a unital
homomorphism $u\colon \ A(D) \to B(\ell_2)$ (necessarily of the form $P\to
P(x)$ for some $x$ in $B(\ell_2)$) such that $\|u\|\le c$ but $\|u\|_{cb}
= \infty$.

The proof of the polynomial boundedness was
simplified in [Kis1] and [DP].

Although this solves the somewhat prototypical case of $A(D)$, it leaves
open the following question: \ is it true that any uniform algebra (\ie 
a unital subalgebra of $C(K)$ for some compact
set $K$) which is proper (\ie $A$ separates
the points of $K$ and  
$A\ne C(K)$) fails (SP)?

See [Kis2] for a partial result on this.
Actually, when $K$ is a domain in $\CC$ with at least 2 holes it is
already unknown in general  whether $\|u\|  =
1$ implies
$\|u\|_{cb}=1$! The case of
a single hole is covered
by [Ag]. See also [DoP] and [Pa4]
for more on this theme.

\n {\bf Remarks.} \item {(i)} See [Ku, KuT] for
recent progress on conditions for an operator
to be similar to a normal operator.

\item {(ii)} The recent paper  [KLM] contains
the following striking example: for any
$c>1$   there is a power bounded operator on
$\ell_2$ which is not similar to any operator 
with powers bounded by $c$. The corresponding
statement for
polynomial boundedness seems open: given
$c>1$, is there a  polynomially bounded 
operator which is not similar to any operator
polynomially bounded by $c$ ?

\medskip
 We now turn to the notion of length
which seems closely connected to the
generalized similarity problem. The ``length''
that we have in mind is analogous to the
following situation: \ consider a unital
semi-group $S$ and a unital generating subset
$B\subset S$, it is natural to say that
$B$ generates $S$ with length $\le d$ if any
$x$ in $S$ can be written as a product $x =
b_1b_2\ldots b_d$ with each $b_i$ in $B$.
We will use a somewhat ``dual'' viewpoint on
the ``length'' based on homomorphisms. Our
main idea can be illustrated in a rather
transparent way on the above simple model of
semi-groups as follows. Assume that $B$
generates $S$ with length $\le d$. Then any
homomorphism $\pi\colon \ S\to B(H)$ (\ie 
$\pi(st) = \pi(s) \pi(t)$ and $\pi(1)=1$)
which is bounded on
$B$ with $\sup\limits_{b\in B} \|\pi(b)\| \le c$ must be bounded on the
whole of $S$ with $\sup\limits_{s\in S} \|\pi(s)\|\le c^d$.

Conversely, assume that we know that for some $\alpha\ge 0$ and $\kappa
\ge 0$, all  homomorphisms $\pi\colon \ S\to B(H)$ satisfy, for some
$c>1$, the following implication:
$$\sup_{b\in B} \|\pi(b)\|\le c \Rightarrow \sup_{s\in S} \|\pi(s)\| \le
\kappa c^\alpha.$$
Then it is rather easy to see that $B$ necessarily generates $S$ with
length $\le [\alpha]$ (integral part of $\alpha$), so that we can replace
$\alpha$ by $[\alpha]$ and $\kappa$ by 1.

We called this a ``dual'' viewpoint because it is reminiscent of the fact
that the closed convex hull $C$ of a subset $B\subset E$ of a Banach space
is characterized by the implication 
$$\sup_{b\in B} f(b) \le 1 \Rightarrow \sup_{s\in C} f(s)\le 1$$
for all continuous real linear forms $f$.

Although this is a wild analogy, we feel that our results on the length
are a kind of ``nonlinear'' analog of this very classical duality
principle for convex hulls.

In [P2], we study various analogs of this concept of length for operator
algebras or even for general Banach algebras. Surprisingly little seems to
have been known up to now. We will now review the main results of our
papers.

\proclaim Definition. An operator algebra $A\subset B({\cal H})$ is said to be of 
length $\le d$ if there is a constant $K$  such that, for any $n$ and any $x$ 
in $M_n(A)$, there is an integer $N = N(n,x)$ and scalar matrices $\alpha_0 \in 
M_{n,N}(\CC)$, $\alpha_1\in M_N(\CC),\ldots, \alpha_{d-1}\in M_N(\CC)$, 
$\alpha_d \in M_{N,n}(\CC)$ together with diagonal matrices $D_1,\ldots, D_d$ in 
$M_N(A)$ satisfying
$$\left\{\eqalign{&x = \alpha_0 D_1 \alpha_1 D_2 \ldots D_d\alpha_d\cr
&\prod^d_0 \|\alpha_i\| \prod^d_1 \|D_i\| \le K\|x\|.}\right.
$$
We denote by $\ell(A)$ the smallest $d$ for which this holds and we call it the 
``length'' of $A$ (so that $A$ has length $\le d$ is indeed the same as $\ell(A) 
\le d$).

Equivalently, we may reformulate this using infinite matrices: \ if we
view as usual $M_n(A) \subset M_{n+1}(A)$ via the mapping $x\to
\left(\matrix{x&0\cr 0&0\cr}\right)$, and if 
we let $\K(A) = \overline{\cup
M_n(A)}$ be the completion of the union with the natural extension of the
norm, then it is easy to check
 that $\ell(A) \le d$ iff any $x$ in $\K(A)$
can be written as
$$x = \alpha_0D_1\alpha_1\ldots D_d\alpha_d$$
with $\alpha_i$ in $\K(\CC)$ and $D_i$
diagonal in $\K(A)$. (The constant
$K$ automatically exists by the open mapping theorem.) 

Our  central result is as follows.

\proclaim Theorem 5 (1999, [P2]). A unital operator algebra $A$ satisfies
(SP) iff $\ell(A) <\infty$. Moreover, let
$$d(A) = \inf\{\alpha\ge 0 \mid \exists~K \ \forall~u \ \|u\|_{cb} \le
K\|u\|^\alpha\}$$
(here of course $u$ denotes an arbitrary unital
homomorphism from $A$ to $B(H)$), then
$$d(A) = \ell(A)$$
and the infimum defining $d(A)$ is attained.

\medskip
\n {\bf Proof that $d(A)\le \ell(A)$.} This is the easy direction. The
converse is much more involved. Assume $\ell(A)\le d$. Consider $x$ in
$M_n(A)$. Recall $\|u\|_{cb} =
 \sup\{\|u_n(x)\|_{M_n(B(H))}\mid
\allowbreak n\ge 1\ 
\|x\|_{M_n(A)}
\le 1\}$. Consider a factorization of the
above form:
$$x = \alpha_0D_1 \ldots D_d\alpha_d$$
with $\alpha_i$ ``scalar'' and $D_i$ ``diagonal''. We have then
$$u_n(x) = \alpha_0 u_N(D_1) \alpha_1\ldots u_N(D_d)\alpha_d$$
hence
$$\|u_n(x)\|\le \prod \|\alpha_i\| \prod 
\| u_N(D_i)\|$$ but clearly since the $D_i$'s
are diagonal $\|u_N(D_i)\| \le \|u\| \
\|D_i\|$ hence
$$\|u_n(x)\| \le \|u\|^d \prod \|\alpha_i\| \prod \|D_i\|$$
which yields (recalling the meaning of $\ell (A)\le d$)
$$\|u\|_{cb}\le K\|u\|^d.$$
\qed\medskip

\n {\bf Remark 6.} Let us briefly return to the derivation problem. If $A$
is a $C^*$-algebra, Kirchberg's argument in [Ki], as slightly improved in
[P2] shows that if we have
 $$\|\delta\|_{cb} \le
\alpha\|\delta\|\leqno(2)$$ for all
$\pi$ and all $\pi$-derivations $\delta\colon \ A\to B(H)$ then we have
$\|u\|_{cb} \le \|u\|^\alpha$ for all $u$ as in Theorem~5. Therefore
$\ell(A)$ is less or equal to the integral part of $\alpha$. This leads us
to conjecture that, in the
$C^*$-case, the best 
possible $\alpha$ in (2) is always an integer.
Also when $A$ is an infinite dimensional
$C^*$-algebra we have no example of $A$ for
which    the best $K$ such that: \
$\forall~u$ \
$\|u\|_{cb}\le K\|u\|^{d(A)}$   is $>1$,
 but we
believe such examples exist (we suspect
$A=B(H)\oplus
\ell_\infty$ might be such an example).

It is easy to see that if $I\subset A$ is a
closed two-sided
 ideal then $\ell(A/I)\le \ell(A)$ and also
that $\ell(A) \le \max\{\ell(I), \ell(A/I)\}.$
If $A$ is a $C^*$-algebra, we have
$$\ell(A) = \max\{\ell(I), \ell(A/I)\}.$$
To show $\ell(I)\le \ell(A)$ we merely use the
fact (due to Arveson) see \eg [Wa] that there
is a ``quasi-central approximate unit" in $I$,
\ie a net $(a_i)$ in the unit ball of $I$
such that for any $x$ in $I$ we have
$xa_i\to  x$ and $a_ix\to  x$ and moreover
(quasi-centrality) $a_ia-aa_i\to  0$ for any
$a$ in $A$.

\n In particular, for all finite sets
$A_1,\ldots, A_n$ of operator algebras we have
$$\ell(A_1\oplus...\oplus A_n) =
\max\{\ell(A_i)\mid 1\le i
\le n\}.$$
The case of infinite direct sums is discussed
in [P6].

\n {\bf Remark 7.} Let $H = \ell_2$. One useful way to apply Theorem~5 is
as follows: \ given a $d$-{\it linear\/} map $w\colon \ A^d\to B(H)$ we
may consider all the possible ways to ``factorize'' $w$ so that there
exist linear bounded maps $v_i\colon \ A\to B(H)$ such that
$$w(a_1,a_2,\ldots, a_d) = v_1(a_1) v_2(a_2)\ldots v_d(a_d).\leqno
\forall(a_1,\ldots, a_d)\in A^d$$
Then we set
$$|||w|||_d = \inf\left\{\prod^d_{i=1} \|v_i\|\right\}$$
where the product runs over all possible ways to ``factorize'' $w$ as
above.  Then let $v\colon \ A\to B(H)$ be a {\it linear\/} map. Assume
that we have a finite set of $d$-linear maps
$w_p$ as before such that
$$v(a_1a_2\ldots a_d) = \sum_p w_p(a_1,\ldots, a_d).\leqno \forall~a_i \in
A\quad (1\le i\le d)$$
Then we set
$$\|v\|_{[d]} = \inf\left\{\sum_p |||w_p|||_d\right\}$$
where the infimum runs over all possible
 ways to write as $v =
\sum\limits_p w_p$. Then if $\ell (A) \le d$,
it is a simple exercise to show that for any
{\it linear\/} $v\colon \ A\to B(H)$ we have
$$\|v\|_{cb} \le K\|v\|_{[d]}.$$
Thus Theorem 5 allows to strengthen the property (SP): \ not only
homomorphisms are c.b.\ but also all linear maps $v$
 for which we have
$\|v\|_{[d]} <\infty$. Actually, it is possible
to show that  $w\to |||w|||_d$ is subadditive
but we will not really  need this.
This norm $\|~~\|_{[d]}$ is closely connected
with the notion of ``multilinear c.b.\ map''
introduced by E.~Christensen and A.~Sinclair
(see [CS1, CS2]).\medskip

\n {\bf Examples.} If $1 <
 {\dim}(A)<\infty$, then $d(A)=1$, so from
now on we assume $\dim(A)=~\infty$. We can now
review the examples of $C^*$-algebras listed
previously:
\item{(i)} If $A$ is commutative $d(A)=2$.
\item{(ii)} If $A = \widetilde \K$ or if $A$
is nuclear, also $d(A)=2$.
\item{(iii)} If $A=B(H)$, then $d(A)=3$.
\item{(iv)} If $A = \widetilde \K\otimes B$
with $B$ arbitrary unital
$C^*$-algebra then $2 \le d(A)\le 3$.
\item{(v)} If $A$ is a $II_1$-factor with property $\Gamma$ then $3 \le
d(A)\le 5$.
\medskip

\n {\bf  Notes:} \ (i) and (ii) are due to
 J.~Bunce and E.~Christensen
(see [C2]). In (iii) $\le 3$ is 
proved in [H1] while $\ge 3$ is proved in
[P2] (see below). (iv) is essentially in [H1].
Finally, concerning (v), Christensen proved in
[C4] that $d(A) \le 44$, but the estimate was
reduced in [P6]. It was also observed in [P6]
that (as pointed out by N.~Ozawa) Anderson's
construction in [An] remains valid on any
$II_1$ factor, thus yielding
$d(A)\ge 3$ for any
$II_1$ factor $A$ by the same argument as in
[P2].

The class of algebras with $d(A) (=\ell(A))$ equal to 2 is closely related
to that of ``amenable Banach algebras'' (see
\eg\ [Pi]). A von~Neumann algebra $M\subset
B({\cal H})$ is called amenable (= injective)
if there is a projection $P\colon \ B({\cal
H})\to M$ with $\|P\|=1$. It is known that a
$C^*$-algebra $A$ is nuclear
($\Leftrightarrow$ amenable by [H2]) iff for
every representation (= $*$-homomorphism)
$\pi\colon \ A\to B(H)$, the von~Neumann
algebra $M_\pi = \pi(A)''$ generated by $\pi$
is amenable (= injective). This motivates the
following

\proclaim Definition. A $C^*$-algebra
 is called semi-nuclear if for any
representation $\pi\colon \ A\to B(H)$ 
generating a semi-finite von Neumann algebra
$\pi(A)''$, the generated algebra $\pi(A)''$
is injective.

\proclaim Theorem 8. ([P2]) For a
$C^*$-algebra $A$, 
$d(A)\le 2$ implies that $A$ is semi-nuclear.

It is an open problem whether in general
semi-nuclear
$\Rightarrow$ nuclear. However, if $A$ is
either the reduced or the full $C^*$-algebra
of a discrete group $G$, then 
$$A \hbox{ nuclear } \Leftrightarrow A \hbox{ semi-nuclear }
\Leftrightarrow G \hbox{ amenable.}$$

  The preceding result  shows that
$d(B(H))>2$, since otherwise $B(H)$ would be
semi-nuclear, which contradicts [An]. Hence, we
have $d(B(H))\ge 3$. Actually, using the length
$\ell(B(H))$ instead, we can obtain a very
simple proof that $d(B(H))= 3$, as follows.

\n {\bf Proof that $\ell(B(H))\le 3$:} \ This very direct proof comes from
[P6]. Fix $n\ge 1$. Let $W_1$ and $W_2$ be any two $n\times n$ unitary
matrices such that
$$|W_1(i,j)| = |W_2(i,j)| = n^{-1/2}.\leqno \forall~i,j$$
Then, for any $x$ in the unit ball of $M_n(B(H))$ (with $\dim H = \infty$)
there are diagonal matrices $D_1,D_2,D_3$ also in the unit ball of
$M_n(B(H))$ such that
$$x = D_1W_1D_2W_2D_3.$$
The proof of this is very simple. Let $S_i$, $i=1,\ldots, n$ be isometries
on $H$ with orthogonal ranges so that
$$S^*_iS_j = \delta_{ij}I.\leqno \forall~i,j$$
Then let
$$D_1(i,i) = S^*_i \quad \hbox{and}\quad D_3(j,j) = S_j$$
and moreover
$$D_2(k,k) = n \sum_{i,j} \overline{W_1(i,k)}S_i x_{ij}S^*_j
\overline{W_2(k,j)}.$$
It is an easy exercise (left to the reader) to check
the announced properties.\qed

By Theorem~1 and Theorem~5, we have:

\proclaim Proposition 9. The Kadison similarity problem has a positive
answer for all unital $C^*$-algebras $A$ iff there is an integer $d_0$ such
that $\ell(A)\le d_0$ for any $C^*$-algebra $A$.

Unfortunately, up to now, the highest known value of $\ell(A)$ for a
$C^*$-algebra is 3, but we conjecture that there are examples of
arbitrarily large length. However, in the {\it non-self-adjoint\/} case, we
have recently been able to prove the following. 

\proclaim Theorem 10. For any integer $d\ge 1$, there is a
(non-self-adjoint) operator algebra $A_d$ such that
$$\ell(A_d) = d.$$

\n {\bf Problem.} Are there uniform algebras with arbitrarily large finite
length?

For uniform algebras no example with $2<\ell(A)<\infty$ is known. However,
it is proved in [P2] that any proper uniform algebra $A$ must satisfy
$\ell(A)>2$. It is also unknown whether there are $Q$-algebras (=
quotients of uniform algebras) $A$ with $2 < \ell(A)<\infty$.\medskip

\n {\bf Sketch of proof of  Theorem 10.} The algebras $A_d$ are not at all
``pathological'', they are the ``obvious'' ones: \ the maximal operator
algebras generated by a sequence of contractions $(x_n)$ to which we
impose the relations
$$x_{n_1}x_{n_2} \ldots x_{n_{d+1}} = 0\leqno ({\cal R}_d)$$
for any $(d+1)$-tuple of integers $(n_1,\ldots, n_{d+1})$. However, while
the proof that $d(A_d)\le d$ is then quite easy, the fact that
$\ell(A_d)>d-1$ has turned out to be much harder to prove. The proof given
in [P4] uses crucially Gaussian random matrices and specifically a recent
difficult estimate due to Haagerup and Thorbj\o rnsen [HT]. We will only
give a brief description
 of the argument from [P4]. Let $P =
P(X_1,X_2,\ldots)$ be a polynomial of degree $\le d$ in {\it
non-commuting\/} (formal) variables $X_1,X_2,\ldots$~. We introduce the
norm 
$$\|P\|_{A_d} = \sup\{\|P(x_1,x_2,\ldots)\|\}\leqno (3)$$
where the supremum runs over all sequences of contractions in $B(\ell_2)$
satisfying $({\cal R}_d)$. It is easy to check that this is a norm of the
set of polynomials $P$ with degree $\le d$. We denote by $A_d$ the
completion of the set of $P$'s equipped with this norm. Clearly, this
defines an operator algebra naturally embedded into $\bigoplus\limits_x
B(H_x)$ where $H_x=\ell_2$ where $x = (x_n)_{n\ge 1}$ runs over the set of
all possible sequences of contractions satisfying $({\cal R}_d)$. In order
to show that $\ell(A_d)>d-1$, the next lemma is crucial. To state it we
first need a specific notation.\medskip

\n {\bf Notation.} Let $H=\ell_2$. Let $m\ge 1$ and $d\ge 1$ be fixed
integers. We will denote by $C(m,d)$ the smallest constant $C$ for which
the following holds: \ if $\{x_i\mid i\in [m]^d\}$ in $B(H)$ satisfies
$$\left\|\sum \lambda_ix_i\right\| \le \sup_{\scriptstyle X_i\in B(H)\atop
\scriptstyle \|X_i\| \le 1} \left\{\left\|\sum \lambda_i X_{i_1}X_{i_2}
\ldots X_{i_d}\right\|\right\} \leqno (4)\quad \forall~\lambda_i \in \CC$$
then $\exists~\hat x_k \in B(H)$, $(1\le k \le m)$ with $\|\hat x_k\| \le
1$ such that
$$x_i  = C\hat x_{i_1} \hat x_{i_2} \ldots \hat
x_{i_d}.\leqno (4)' \quad
\forall~i\in [m]^d$$

\proclaim Lemma 11. For any $m\ge 1$ and $d\ge 1$, we have
$$\delta_d m^{d-1\over 2} \le C(m,d) \le m^{d-1\over 2}$$
where $\delta_d>0$ is a constant independent of $m$.

\n {\bf Example.} In the case $d=2$, this means
 the following: \ if
$x_{ij} \in B(H)$ $(i,j=1,...,m)$ satisfy
$$\left\|\sum_{ij\le m}
\lambda_{ij}x_{ij}\right\| \le \sup_{\|X_i\|\le
1}
\left\|\sum \lambda_{ij}X_iX_j\right\|\leqno \forall~\lambda_{ij} \in \CC$$
then $x_{ij}$ can be factorized as
$$x_{ij} = C\hat x_i \hat x_j \quad \hbox{with}\quad \|\hat x_i\|\le
1\leqno (5)$$
but in general the best possible $C$ will be $\simeq \sqrt m$. This case
is rather easy to prove given the state of the art. However, already the
case $d=3$ is more delicate, and as we already mentioned the case of an
arbitrary $d$ requires the upper estimates given in [HT] which are highly
non-trivial. An easier proof of the
lower bound (which is the difficult part) in
Lemma~11 would be most welcome. 
\medskip

\n {\bf Remark 12.} Given $\{x_i\mid i\in [m]^d\}$ in $B(H)$ satisfying
(4), we can define a {\it linear\/} map
$$v\colon \ A_d\to B(H)$$
by setting $v(X_{i_1}X_{i_2}\ldots X_{i_d}) = x_{i_1i_2\ldots i_d}$ with
$1\le i_1,i_2,\ldots, i_d\le m$ and $v(X_{i_1}X_{i_2} \ldots X_{i_k}) =~0$
in all other cases.
Then it can be shown, using the
factorization
of multilinear cb maps of
Christensen-Sinclair and Paulsen-Smith (see
[P4]) that
$\|v\|_{cb}$ is equal to the smallest constant
$C$ such that (4)$'$ holds.

We now wish to sketch how Lemma~11 is used to prove that $\ell(A_d) >d-1$.
To lighten the exposition, we will restrict to
the simplest case:
\
$d=3$. So we will show that Lemma~11 implies
$\ell(A_3)>2$. We will show that if
$\ell(A_3)\le 2$ then $C(m,d)\le K\sqrt m$ for some $K$, but this will
contradict Lemma~11 for $d=3$ since $(d-1)/2 = 1 > 1/2$, whence the
conclusion that $\ell(A_3)>2$.

Now assume $\ell(A_3)\le 2$. Let $\{x_{i_1i_2i_3} \mid i\in [m]^3\}$ be as
in the definition of $C(m,d)$ for $d=3$. For convenience, we extend the
function $(i_1,i_2,i_3)\to x_{i_1i_2i_3}$ to
the whole of $\NN^3$ by setting it equal to
zero outside $[1,\ldots, m]^3$. 

We will use Remark 7.

Let $v\colon \ A_3\to B(H)$ be the linear map defined by $v(1) = 0$,
$v(X_i)=0$, $v(X_{i_1}X_{i_2}) = 0$ and finally:
$$v(X_{i_1}X_{i_2}X_{i_3}) = x_{i_1i_2i_3}.$$
It is easy to see using (3) and (4) that
$$\|v\|\le 1.$$
We claim that (4) implies (with the notation
of Remark 7)
$$\|v\|_{[2]} \le 2+2\sqrt m.$$
We will use the following notation: \ we
consider the disjoint union
$$\Omega  = {\phi}\cup \NN \cup \NN^2 \cup
\NN^3,$$ and we set
$$\eqalign{X^\phi &= 1\cr
X^i &= X_i \quad \hbox{if}\quad i\in \NN\cr
X^{ij} &= X_iX_j \quad \hbox{if}\quad (ij) \in \NN^2\cr
X^{ijk} &= X_iX_jX_k \quad \hbox{if}\quad (ijk)\in \NN^3.}$$
For $i\in\Omega$ we set $|i| = 0$ if $i=\phi$,
and $|i|=k$ if $i\in
\NN^k$.

\n With this notation any polynomial $P$ in
$A_3$ can be written as a finite sum
$$P = \sum_{i\in\Omega} \lambda_i(P)X^i$$
with $\lambda_i(P)\in \CC$. We have then 
$v(X^i)=x_i$ for all $i$ in
$\Omega$, hence $\forall~P_1,P_2\in A_3$
$$v(P_1P_2) = \sum_{i,j\in\Omega} \lambda_i(P_1) \lambda_j(P_2) x_{ij}$$
where $ij$ denotes now the multi-index of length $\le 6$ obtained by
putting $j$ after $i$. We set $|ij| = |i| + |j|$. But since $x_{ij} = 0$
unless $|i| + |j|=3$ we find a decomposition of $v$ as follows:
$$v(P_1P_2) = \sum_{(\alpha\beta)\in J} w_{\alpha\beta} (P_1,P_2) \leqno
(6)$$
where the sum runs over the set $J$ of all pairs $(\alpha\beta)$ in
[0,1,2,3] such that $\alpha+\beta=3$, and where  $w_{\alpha\beta}$ are
bilinear forms on $A_3\times A_3$ defined by setting:
$$w_{\alpha\beta}(P_1,P_2) = \sum_{\scriptstyle |i|=\alpha \atop
\scriptstyle |j|=\beta} \lambda_i(P_1) \lambda_j(P_2) x_{ij}.$$
Using (4) it is easy to see that if $(\alpha\beta)$ is either 
$(3~ 0)$ or $(0~ 3)$
then with the notation of Remark 7
$$|||w_{\alpha\beta}|||_2\le 1.$$
The remaining possibilities in $J$ are only
(2~1) and (1~2). But if
$(\alpha\beta) = (1~~2)$ for instance we can write
$${w_{\alpha\beta}(P_1,P_2)  = 
\left(\sum_{|i|=1} \lambda_i(P_1) e_{1i}
\otimes I\right)
\phantom{=} \left(\sum^m_{k=1} e_{k1} \otimes
\sum_{|j|=2} \lambda_j(P_2) x_{kj}\right)}$$
(here we identify $B(H)$ with $B(H) \overline\otimes B(H)$ and denote by
$(e_{ki})$ the standard matrix units in
$B(H)$). Using this, one
can check rather easily that if
$(\alpha\beta) = (2~1)$ or $(1~2)$ then
$$|||w_{\alpha\beta}|||_2 \le \sqrt m.$$
Thus using (6) we obtain our claim that
$$\|v\|_{[2]} \le 2+2\sqrt m.$$
Then if we assume $\ell(A_3) \le 2$, Remark~7 ensures that
$$\|v\|_{cb} \le K(2+2\sqrt m).$$
Now by Remark~12 this implies that $\{x_i\mid i\in [m]^3\}$ satisfies
(4)$'$ with $C\le K(2+\sqrt m)$. Thus we
conclude that $C(m,3) \le K(2+\sqrt m)$ but
this obviously contradicts Lemma~11 with
$d=3$. Thus we have shown, by this
contradiction, that $\ell(A_3)>2$.
\qed
\bigskip

\magnification\magstep1
\baselineskip = 18pt
\overfullrule = 0pt

\def\n{\noindent}
\def\qed{\hfill {\rm qed}}

\def\calk{{\cal K}}

\def\CC{ \;Ê{}^{ {}_\vert }\!\!\!{\rm C}}

\def\NN{{\rm I}\!{\rm N}}

\def\K{\calk}

The notion of length is quite natural in the more general context of a Banach
algebra $B$ generated by a family of subalgebras $B_i\subset B$ $(i\in I)$.
For simplicity, we will restrict ourselves to the case of a pair of
subalgebras $B_1\subset B$, $B_2\subset B$. In this case, we say that
$B_1,B_2$ generate $B$ with length $\le d$ if
 there is a bounded subset of the union  $C
\subset B_1\cup B_2$ such that every $x$ in the unit ball of $B$ belongs to
the closed convex hull of   the union $\cup_{j=1}^d C^j$ where
$$C^j = \{x_1x_2\ldots x_j\mid x_k\in C \quad \forall k=1,...,j
.\}$$
  Assuming that this holds, let
$u\colon \ B\mapsto \beta$ be a continuous homomorphism into another Banach
algebra $\beta$. It is then easy to check that 
$$\|u\| \le K \sum^d_{j=1} \max\{\|u_{|B_1}\|,
\|u_{|B_2}\|\}^j$$ where $K$ is a constant (depending only on
$d$ and the size of the subset $C$).

\n Thus if $B,\beta$ and $u$ are all unital we obtain (since all
the norms are now $\ge1$)
$$\|u\| \le dK \max\{\|u_{|B_1}\|, \|u_{|B_2}\|\}^d.$$

In the converse direction, assuming again
$B_1,B_2$ and $B$ all unital, let $\hbox{alg}(B_1,B_2)$ denote the algebra
generated by $B_1$ and $B_2$, which we assume is dense in $B$. Assume that
every unital homomorphism $u\colon \ \hbox{alg}(B_1,B_2)\to \beta$ into an
arbitrary unital Banach algebra $\beta$ such that $\|u_{|B_1}\|<\infty$ and
$\|u_{|B_2}\|<\infty$ is actually bounded and satisfies
$$\|u\| \le K(\max\{\|u_{|B_1}\|, \|u_{|B_2}\|\})^\alpha$$
where $K$ and $\alpha\ge 0$ are independent of $u$ and $\beta$.

\n Then it follows (see [P2, \S 8]) that $B_1,B_2$ generate $B$
with length at most equal to the integral part of $\alpha$. For
example, let $A$ be a unital operator algebra, and let $B =
\K(A)$. We may consider the subalgebra
$B_1\subset B$ formed of all the diagonal matrices (viewing the
elements of $\K(A)$ as bi-infinite matrices with coefficients in
$A$) and we  let $B_2 = \K(\CC)$. 

\n It is then easy  to check
that $\ell(A) \le d$ implies that $B_1,B_2$ generate $B$ with
length $\le 2d+1$. Conversely, if $B_1,B_2$ generate $B$ with
length $\le m$, then $\ell(A) \le \left[{m+1\over 2}\right]$. 
\medskip

\n {\bf Remark.} The slight discrepancy appearing here comes from the fact
that in the products appearing in the subset $C^d$ we do not specify that the
first term of the product must lie in $B_2$ or $B_1$ while in the
corresponding definition of $\ell(A)$ the analogous term must be in $B_2$.
This difficulty can be circumvented: \ one should then consider homomorphisms
$u\colon \ \hbox{alg}(B_1,B_2) \to \beta$ such that $\|u_{|B_2}\| = 1$ and
study the inequality $\|u\| \le K\|u_{|B_1}\|^\alpha$. See [P3]
for more variations on this theme.

The case
study of $\ell(A)$ suggests to examine many other examples of the same kind,
for instance the pair $B_1 = \K(A_1)$, $B_2 = \K(A_2)$ where
$A_1
\subset A$,
$A_2\subset A$ are two closed subalgebras. In particular, we may consider the
case where $A$ is the maximal tensor product of two unital $C^*$-algebras
$C_1,C_2$: namely we take
$A= C_1\otimes_{\max} C_2$ with $A_1=C_1\otimes 1$ and $A_2 =
1\otimes C_2$. All these cases are studied in [P3], to which we
refer the reader for several illustrating examples and more
information.

 \bigskip\bigskip

\centerline {\bf References}

\item{[Ag]} J. Agler. Rational dilation on an
annulus. Ann. of Math.   121 (1985) 
537--563.

\item{[An]}
  J. Anderson. Extreme points in sets of positive linear maps on
${\cal B}({\cal H})$. J. Funct. Anal. 31 (1979) 195-217.

\item{[BP]}  D. Blecher and V. Paulsen. 
Explicit construction of universal operator
algebras and applications to polynomial
factorization. Proc. Amer. Math. Soc.   112
(1991) 839-850.

\item{[C1]}   E. Christensen.
    Extensions of derivations.   	 J.
Funct. Anal.   27  (1978)  234--247 

\item{[C2]} $\underline{\hskip1.5in}$. 
Extensions of derivations II   Math.
Scand.  50 (1982)    111--122.

\item{[C3]}  $\underline{\hskip1.5in}$. On non
self adjoint representations of operator
algebras  Amer. J. Math. 103 (1981) 817-834.

\item{[C4]}   $\underline{\hskip1.5in}$.
Similarities of $II_{1}$ factors with property
$\Gamma$.  Journal Operator Theory  15 (1986)
281-288. 

\item{[C5]}  
$\underline{\hskip1.5in}$. Perturbation of
operator algebras II. Indiana Math.\ J. {  26}
(1977), 891--904.

\item{[CS1]} E. Christensen and A. Sinclair. 
Representations of completely bounded multilinear operators.
J. Funct. Anal. 72 (1987) 151-181.

\item{[CS2]}   $\underline{\hskip1.5in}$.  A
survey of completely bounded operators.  Bull.
London Math. Soc.  21 (1989) 417-448.

\item{[DP]}  K. Davidson and V. Paulsen.
On polynomially bounded operators,
J. f\"ur die reine und angewandte Math. {  487}
(1997) 153-170.

\item{[DoP]} R. Douglas   and V. Paulsen.
{\it Hilbert modules over function algebras. 
}  Pitman Longman 1989.

\item{[DW]} H. G. Dales and  W. H. Woodin. {\it
An
introduction to independence for analysts.}
  London Mathematical
Society Lecture Note Series, 115. Cambridge
University Press, Cambridge-New York, 1987. 

\item{[EM]}
 L. Ehrenpreis  and F.I. Mautner.  Uniformly
bounded representations of groups  Proc. Nat.
Acad. Sc. U.S.A. 41 (1955) 231-233.

\item{[H1]}  U.  Haagerup.  Solution of the
similarity problem for cyclic representations
of $C^*$-algebras.  Annals of Math.  118
(1983), 215-240.

\item{[H2]} $\underline{\hskip1.5in}$. 
All  nuclear $C^*$-algebras are amenable.
 Invent. Math.   74 (1983) ) 305--319.

\item{[Ha]} P. Halmos. Ten problems in Hilbert
space. Bull Amer. Math. Soc.  (1970) 

 \item{[HT]} U. Haagerup and S. Thorbj{\o}rnsen.
Random matrices and $K$-theory for exact
$C^*$-algebras. Documenta Math. 4 (1999)
341-450.

\item{[Ka]} R.  Kadison  On the orthogonalization
of operator representations. Amer. J. Math.  77
(1955) 600-620.

   \item{[Ki]} E. Kirchberg.  The derivation and the similarity problem 
are equivalent.
 J. Operator Th. 36 (1996) 59-62.

\item{[Kis1]} S. Kislyakov. Operators that are
(dis)similar to a contraction: Pisier's
counterexample in terms of singular integrals.
Zap. Nauchn. Semin. S.- Peterburg. Otdel. Mat.
Inst. Steklov (POMI).

   \item{[Kis2]} $\underline{\hskip1.5in}$. 
Similarity problem for certain martingale
uniform algebras. Preprint. 
 
 \item{[KLM]} N. Kalton and C. Le Merdy.
Solution of a problem of Peller concerning
similarity. J. Op. Theory. To appear.

 \item{[Ku]} S. Kupin. Similarit\'e \`a un
op\'erateur normal et certains probl\`emes
d'interpolation. Th\`ese, Universit\'e de
Bordeaux I, (3 avril 2000).

\item{[KuT]} S. Kupin and S. Treil. Linear
resolvent growth of a weak contraction does not
imply its similarity to a normal operator.
Preprint 2000, to appear. 

\item{[Pa1]} V. Paulsen.   Completely bounded maps and
dilations.  Pitman Research Notes in
Math. 146, Longman, Wiley, New York, 1986.

 \item{[Pa2]}  $\underline{\hskip1.5in}$. 
Completely bounded homomorphisms of operator
algebras.  Proc. Amer. Math. Soc. 92 (1984)
225-228.

 \item{[Pa3]} $\underline{\hskip1.5in}$. 
Every completely  polynomially bounded
operator is similar to a contraction.  J.
Funct. Anal. 55 (1984) 1-17. 

 \item{[Pa4]}  $\underline{\hskip1.5in}$.
 Toward a theory of $K$-spectral sets. Surveys
of some recent results in operator theory,
Vol. I, 221--240, Pitman Res. Notes Math.
Ser., 171, Longman Sci. Tech., Harlow, 1988.

\item{[P1]}  G. Pisier.  A polynomially bounded
operator on Hilbert space which is not similar
to a contraction. Journal Amer. Math. Soc. 10
(1997) 351-369. \vskip 12pt

\item{[P2]}$\underline{\hskip1.5in}$.  The
similarity degree of an operator algebra.
   St. Petersburg Math. J. 10 (1999) 103-146.
 
\item{[P3]}   $\underline{\hskip1.5in}$. Joint
similarity problems and the generation of
operator algebras with bounded length. Integr.
Equ. Op. Th.
 31 (1998) 353-370.

\item{[P4]}  $\underline{\hskip1.5in}$. The
similarity degree of an operator algebra II.
 Math. Zeit. (2000)  To appear.

\item{[P5]}   $\underline{\hskip1.5in}$.
Similarity problems and completely bounded maps.
Springer Lecture notes 1618 (1995).
 
\item{[P6]}   $\underline{\hskip1.5in}$. 
Remarks on the similarity degree of an
operator algebra. Preprint.
 
\item{[P7]} $\underline{\hskip1.5in}$. Are
unitarizable groups amenable ? preprint, 1999.

 \item{[P8]} $\underline{\hskip1.5in}$.
 Operator spaces and similarity problems,
Proc. Berlin ICM  1998, Doc. Math. Extra Vol.
I, 429-452.

   \item{[Pi]}  J.P. Pier. {\it   Amenable
Banach algebras.} Pitman, Longman 1988. 

 \item{[SS]} A. M. Sinclair and  R.R. Smith. 
 {\it Hochshild cohomology of von Neumann
algebras.} L. M. S. Lecture Notes series.
 Cambridge University Press, Cambridge 1995.

  \item{[Wa]} S. Wassermann. {\it Exact
$C^*$-algebras and related topics.}
Lecture Notes Series, Seoul Nat. Univ. (1994).
 
\vskip12pt

Texas A\&M  University

College Station, TX 77843, U. S. A.

and

Universit\'e Paris VI

Equipe d'Analyse, Case 186,
 
75252 Paris Cedex 05, France
\end

\bye